\documentclass[letterpaper, 10 pt, conference]{ieeeconf} 
\IEEEoverridecommandlockouts
\usepackage{enumerate}
\usepackage{cite}
\usepackage{amsmath,amssymb,amsfonts}
\usepackage{graphicx}
\usepackage{textcomp}
\usepackage{fancyhdr}
\usepackage{array,booktabs}
\usepackage{xcolor}
\usepackage{siunitx}
\usepackage[short]{optidef}
\usepackage[
pdfauthor={Jakob Harzer},
pdftitle={Integration Error Regularization in Direct Optimal Control using Embedded Runge Kutta Methods},
pdfstartview=XYZ,
bookmarks=true,
colorlinks=true,
linkcolor=black,
urlcolor=blue,
citecolor=blue,
pdftex,
bookmarks=true,
linktocpage=true, 
hyperindex=true,
]{hyperref}

\newcommand{\R}{{\mathbb{R}}}
\newcommand{\order}[1]{\mathcal{O}(#1)}

\newcommand{\jochem}[1]{{\color{blue}#1}}

\begin{document}
\title{Integration Error Regularization in Direct Optimal Control using Embedded Runge Kutta Methods
\thanks{This research was supported by DFG via Research Unit FOR 2401, project 424107692, 504452366 (SPP 2364), and 525018088, by BMWK via 03EI4057A and 03EN3054B, and by the EU via ELO-X 953348.
\{\textbf{jakob.harzer}, jochem.de.schutter, moritz.diehl\}@imtek.uni-freiburg.de.}}

\author{Jakob Harzer, Jochem De Schutter, Moritz Diehl}

\maketitle
\thispagestyle{empty} 
\begin{abstract}
    In order to solve continuous-time optimal control problems, direct methods transcribe the infinite-dimensional problem to a nonlinear program (NLP) using numerical integration methods.
    In cases where the integration error can be manipulated by the chosen control trajectory, the transcription might produce spurious local NLP solutions as a by-product.
    While often this issue can be addressed by increasing the accuracy of the integration method, this is not always computationally acceptable, e.g., in the case of  embedded optimization.
    Therefore, alternatively, we propose to estimate the integration error using established embedded Runge-Kutta methods and to regularize this estimate in the NLP cost function, using generalized norms.
    While this regularization is effective at eliminating spurious solutions, it inherently comes with a loss of optimality of valid solutions.
    The regularization can be tuned to minimize this loss, using a single parameter that can be intuitively interpreted as the maximum allowable estimated local integration error.
    In a numerical example based on a system with stiff dynamics, we show how this methodology enables the use of a computationally cheap explicit integration method, achieving a speedup of a factor of 3 compared to an otherwise more suitable implicit method, with a loss of optimality of only 3\%.

\end{abstract}

\section{Introduction}

In this paper, we consider the use of direct methods to solve optimal 
control problems (OCPs) of the form 
\begin{mini!}
    {x(\cdot), u(\cdot)}{\int_{0}^{t_\mathrm{f}} l(x(t),u(t)) \mathrm{d}t + E(x(t_{\mathrm{f}}))\label{eq:general_OCP_cost}}
    {\label{ocp:general_OCP}}{}
    \addConstraint{0}{= x(0) - \bar x_0}\label{eq:general_OCP_initial}
    \addConstraint{\dot x(t)}{= f(x(t), u(t)), \qquad&&\forall t\in [0,t_\mathrm{f}]}\label{eq:general_OCP_dynamics}
    \addConstraint{0}{\leq h(x(t), u(t)), \qquad&&\forall t\in [0,t_\mathrm{f}]}\label{eq:general_OCP_constr}
\end{mini!}
with independent variable $t \in [0, t_\mathrm{f}]$, state $x(t) \in \R^{n_x}$, control $u(t) \in \R^{n_u}$, cost \eqref{eq:general_OCP_cost}, state dynamics \eqref{eq:general_OCP_dynamics}, state and control constraints \eqref{eq:general_OCP_constr}, and initial state $\bar{x}_0$.
This involves first discretizing the continuous trajectories for state and control, and then solving the resulting nonlinear programming problem (NLP).
For this, the solution maps of the dynamics are approximated using numerical integration schemes whose numerical accuracy possibly depends on the control trajectory.
In this case, a side effect of the discretization can be the existence of spurious solutions, characterized by large integration errors.
We quote \cite{Betts2010}: ``If there is a flaw in the problem formulation, the optimization algorithm will find it.''



\medskip\noindent\textbf{(Minimal Example)} Consider the optimal control problem
\begin{mini!}%
    {x(\cdot), u \in \R}{1 - x(1)}%
    {\label{eq:simpleExample_OG}}{}%
    \addConstraint{x(0)}{=\bar x_0}%
    \addConstraint{\dot{x}(t)}{= -u x, \qquad \forall t \in [0,1]}%
    \addConstraint{0}{\leq u \leq 30}%
\end{mini!}%
with the state $x\in \R$, the initial value $\bar{x}_0 = 1$, and a constant control $u \in \R$.
The optimal solution of \eqref{eq:simpleExample_OG} is trivially given by $u^* = 0$.
Using a direct approach to solve this problem, we discretize the continuous state trajectory $x(\cdot)$ using a single step of a Gauss-Legendre collocation scheme with 4 stages (GL8) and solve the resulting nonlinear program (NLP). 
For two different initializations for the variable $u$,  i.e., $u_{0,a} = 0$ and $u_{0,b} = 10$, we obtain two different solutions: $u^*_a = 0$, and $u^*_b = 30$.

\begin{figure}
    \includegraphics[width=1\linewidth,trim=0.53cm 0.4cm 0 0.2cm, clip]{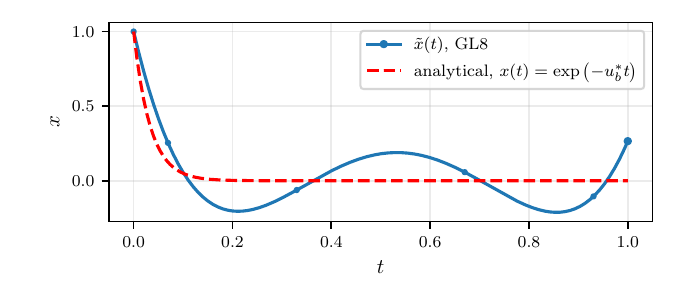}
    \caption{Spurious optimal state trajectory of the discretized problem \eqref{eq:simpleExample_OG} using a 4-stage Legendre collocation integration scheme, and with initialization $u_{0,b} = 10$.}\label{fig:SimpleExample_LocalMinima}
\end{figure}

Fig. \ref{fig:SimpleExample_LocalMinima} shows the optimal state trajectory of the latter solution.
Here, the optimizer converged to a spurious minimum, with large $u$ and hence, very stiff dynamics and a large integration error. 
\begin{figure}
    \centering
    \includegraphics[width=1\linewidth,trim=0.2cm 0.4cm 0 0.3cm, clip]{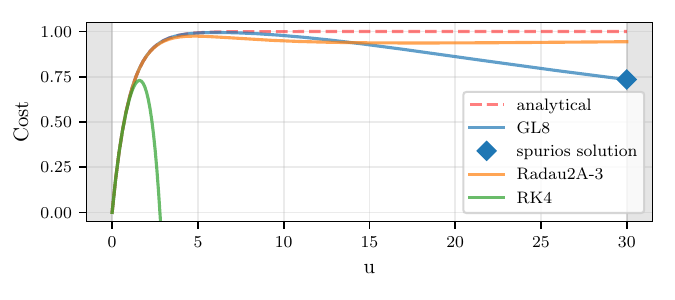}
    \caption{NLP cost function over the feasible set for different integration methods}
    \label{fig:SimpleExample_Objective}
\end{figure}
Fig.~\ref{fig:SimpleExample_Objective} illustrates how the NLP cost function produced using the GL8 integration method allows for an artificial local minimum that is not present in the analytical cost function. The plot shows how the same holds for different integration methods as well,  such as, e.g., RadauIIa with 3 stages (Radau2A-3) or Runge-Kutta 4 (RK4).



\subsection{Related Work}
Typical remedies to address this problem are: (a) add \textit{constraints or regularization} to the optimization variables, (b) provide a \textit{better initialization} to the variables, (c) most obviously, choose a \textit{more accurate discretization} for the given problem or (d) add a \textit{mesh refinement mechanism}. 
These methods iteratively update the discretization mesh for state and control, until a specified solution accuracy is reached.
After each solution of the NLP, the integration error is estimated and the mesh updated \cite{Betts1998, Betts2010, Darby2011a,Darby2011b}. 
Here, error estimates are typically obtained by comparison of interval endpoints with a higher-accuracy simulations \cite{Betts1998}, or by approximating the integral of the absolute violation of dynamics of the continuous representation of the solution \cite{Betts2010,Darby2011a}.

In a scenario where the integration error is large and a more accurate discretization or mesh refinement not feasible due to a limited computational budget, for example in an online setting such as model predictive control, approaches such as a non-uniform grid \cite{Quirynen2015} can be helpful.
An approach similar to the one proposed in this paper is taken by \cite{Lahr2024}, where instead of deterministic integrators, probabilistic integrators are employed and their estimated uncertainty is penalized.
Also quite similar is the ``Integral-Penalty-Barrier Method'' \cite{Neuenhofen2020, neuenhofen2021} where the integral of the squared dynamics violation is added to the cost function instead of solving the collocation constraints exactly.

\subsection{Contribution} 

While the above remedies (a) and (b) try to deal with the symptoms of the problem, we present another simple approach that mitigates the root cause of the problem: estimate and penalize the integration error online. 
The premise is that the regularization of the estimated integration error leads to a control trajectory for which the solution of the NLP is of higher accuracy.
To estimate the error, we employ embedded Runge-Kutta methods, and modify the cost function of the original problem using generalized norms of the estimate.
The approach can be used to stabilize the solution of otherwise unsuited discretization schemes and to trade off optimality and accuracy.

\subsection{Structure}
This paper is structured as follows. Section 2 summarizes some preliminaries on direct optimal control and error estimates. Section 3 presents the proposed regularization method, Section 4 showcases the capability with a numerical example, and Section 5 provides a conclusion.


\section{Discretization and Background}

We discretize problem \eqref{ocp:general_OCP} into an NLP using the well established direct multiple shooting methodology \cite{Rawlings2017, Betts2010, Biegler2010}. 
We divide state and control trajectories into $N$ equidistant intervals of size $h = t_\mathrm{f}/N$ and introduce shooting nodes $x_0,...\,,x_N$ that approximate the state trajectory at the interval boundaries as $x_k \approx x(t_k), \forall k \in \mathcal{K} = \{0,1, \dots, N-1 \}$, with $t_0 = 0 \leq t_1 \leq \dots \leq t_N = t_\mathrm{f}$,
and consider a piecewise constant control parametrization as $u(t) = u_k$ for $t \in [t_k,t_{k+1}),k \in \mathcal{K}$.

\subsection{Runge-Kutta Methods}

To numerically approximate the local solution of the state dynamics $x'(t), t \in [t_k, t_{k+1}]$  with the initial point $x'(t_k) = x_k$ on an interval $k$, we consider a single integration step of an implicit Runge-Kutta (IRK) method of size $h$. 
An IRK method with $d$ stages, with stage values $z_{k,i} \in \R^{n_x}$ and dynamics evaluations $k_i = f(z_{k,i}, u_k)$, is defined by the system of equations
\begin{subequations}
    \begin{align}
        z_{k,i} &= x_k + h \sum_{j=1}^{d} a_{i,j} k_j  \\
        x_{k+1} &= x_k + h \sum_{i=1}^{d} b_i k_i 
    \end{align}
\end{subequations}
which we summarize in the implicit expression
\begin{align}
    0 =  G_\mathrm{IRK}(x_k, x_{k+1}, u_k, z_k) \ ,
\end{align}
with $z_k = \begin{bmatrix} z_{k,1}, \ldots, z_{k,d}\end{bmatrix}^{\top}$.
The coefficients of the IRK scheme are usually summarized in a Butcher tableau which consists of the vectors $c \in \R^d$  and $b \in \R^d$, as well as the matrix $A \in \R^{d \times d}$. 
These coefficients are chosen such that the IRK method is of order $p$.

The local error of a single interval is given by
\begin{align}
    x_{k+1} - x'(t_k + h)  =  \order{h^{p+1}}
\end{align}
which then accumulates over the horizon to a global error
\begin{align}
    x_{N} - x(t_\mathrm{f}) =  \order{h^{p}}
\end{align}
where here $x(\cdot)$ is the solution to the initial value problem \eqref{eq:general_OCP_initial}-\eqref{eq:general_OCP_dynamics} for a given control trajectory.
For more details on IRK methods, we refer the reader to \cite{Hairer1996,Hairer2008}
\subsection{NLP Formulation}
The NLP is the given by
\begin{mini!}
    {w}{J(w)}
    {\label{eq:NLP_orginal}}{}
    \addConstraint{x_0}{=  \bar x_0}
    \addConstraint{0}{      = G_\mathrm{IRK}(x_k, x_{k+1}, u_k, z_k), \quad & \forall k \in \mathcal{K} }
    \addConstraint{0}{\leq h(x_k, u_k, z_k),  \quad & \forall k \in \mathcal{K}}
\end{mini!}
where the variables $w$ consist of the shooting nodes, the integration variables and the control variables:
\begin{align}
    w = (x_0,...\,,x_N, z_0, ...\,, z_{N-1}, u_0, ...\,, u_{N-1}) \in \R^{n_w}.
\end{align}
The cost function $J: \R^{n_w} \rightarrow \R$ is obtained using a numerical approximation of the integral in Equation \eqref{eq:general_OCP_cost}, for example by quadrature, using the same integration method as for the dynamics.
For more information on direct methods for optimal control problems, we refer to \cite{Betts2010,Biegler2010,Rawlings2017}.

\subsection{Error Estimate via Embedded Runge Kutta Methods}

A widely used error estimate of the integration error compares the endpoint $x_{k+1}$ of the integration methods to a second integrator result $\hat{x}_{k+1}$ of typically lower order $\hat{p} \leq p$.
For an efficient computation, the `embedded' integrator uses a different set of weights $\hat{b}$ but the same stage dynamics evaluations $k_i$:%
\begin{align}%
    \hat{x}_{k+1} = x_k + h \sum_{i = 1}^{d} \hat{b}_i k_i %
\end{align}%
The difference between the two endpoints is then an estimate for the local integration error (of the less accurate embedded method \cite{Hairer2008}) in the $k$-th interval:
\begin{subequations}%
    \begin{align}%
        \hat{e}_k(x_k, u_k, z_k) &= \hat{x}_{k+1} - x_{k+1}\\
            &= \sum_{i = 1}^{d} (\hat{b}_i - b_i) k_i
    \end{align}\label{eq:embIntErrApprox}%
\end{subequations}%
The result of the higher-order integration scheme $x_{k+1}$ is used to `continue the integration' \cite{Hairer2008}.
Runge-Kutta integration methods that also provide a second integration result are known as `embedded integrator pairs' of order $p(\hat{p})$.
The embedded integrator is usually defined by extending the Butcher tableau by an additional row, as Fig. \ref{fig:EmbeddedButcherTableaus} shows.
\begin{figure}
    \renewcommand{\arraystretch}{1.2}
    \begin{align*}
        \begin{array}{c|c}
            c & A \\
        \hline
            & b^\top \\ & \hat{b}^\top
    \end{array} &&
        \begin{array}{c|c c}
            0 & 0 & 0 \\
            1 & 1 & 0  \\ \hline
              & \frac{1}{2} & \frac{1}{2} \\ 
              & 1 & 0
        \end{array} &&
        \begin{array}{c|cc}
            c_0 & 0 & 0 \\
            c & 0 & A \\
            \hline
               &0 & b^\top \\ 
               &\gamma_0 & \hat{{b}}^\top
        \end{array}
    \end{align*}
    \caption{Left: Butcher tableau of an embedded integrator pair. Middle: Heun-Euler of Order 2(1). Right: extension of an existing IRK method with an embedded integrator.}\label{fig:EmbeddedButcherTableaus}
\end{figure}
Embedded integrators are primarily employed for simulations with adaptive step size control. 

The most simplest explicit embedded integrator pair is the Heun-Euler method which combines the trapezoidal (Heun) scheme of order 2 and an explicit Euler of order 1 for error estimation (middle table in Figure \ref{fig:EmbeddedButcherTableaus}).
Also well known are the Fehlberg method of order 4(5) with 6 stages, (not to be confused with the popular RK4 method), and the Dormand-Prince method of order 5(4) with 7 stages, used in the popular \texttt{ode45} solver.
Their Butcher tableaus can be found in \cite{Hairer1996,Hairer2008}.

As shown by \cite{Hairer1996, Kouya2013}, general $d$-stage \textit{implicit} RK methods with coefficients $c,A,b$ (whose matrix $A$ has at least one real eigenvalue) can be extended with a lower order ($\hat{p} = d$) embedded method, by adding an additional stage $c_0 = 0$, as shown in the right table in Figure \ref{fig:EmbeddedButcherTableaus}.
Here, $\gamma_0$ is any nonzero constant and the coefficients $\hat{b}$ can be found by solving a linear system of equations.
We refer to \cite{Kouya2013} for details on how to compute these coefficients.

\section{Embedded Residual Minimization}

Let $\hat{e}_k \in \R^{n_x}$ be an estimate for the local integration error in integration interval $k$ that whose component-wise absolute value we would like to keep bounded by a given maximum error $e_{k,\mathrm{max}}\in \R^{n_x}$.
In this paper, we will focus on the estimate \eqref{eq:embIntErrApprox} obtained using an embedded Runge-Kutta method, which is an approximation of the local integration error of the lower-order embedded integrator.
We collect the error terms in the vectors
\begin{align}
    \hat{E}(w) = [\hat{e}_0,\dots,\hat{e}_{N-1}]^\top \in \R^{N n_x}\label{eq:E_hat}
\end{align}
and similarly $ E_\mathrm{max} = [e_{0,\mathrm{max}},\dots,e_{N-1,\mathrm{max}}]^\top \in \R^{N n_x}$.
We propose to normalize the error $\hat{e}_k$ over each interval using a weight vector $\omega \in \R^{n_x}$, such that the only single tuning parameter is the relative maximum error $e_{\mathrm{max}} \in \R$, such that we have
\begin{align}
    E_\mathrm{max}(e_{\mathrm{max}}) = e_{\mathrm{max}}\cdot[\omega, \dots, \omega]^\top\label{eq:Emax_singleTuning}
\end{align}
The error estimate \eqref{eq:E_hat} is then used in the additional regularization cost function $\phi: \R^{n_w} \rightarrow \R$\jochem{:}
\begin{align}
    \phi(w; E_\mathrm{max}, p, q) = \|\mathrm{diag}(E_\mathrm{max})^{-1}\hat{E}(w)\|^q_p\label{eq:pqregularization}
\end{align}
where we typically we choose $p,q = 2$.
We use it to modify the cost function of the NLP \eqref{eq:NLP_orginal}, which now reads as:%
\begin{mini!}%
    {w}{J(w) + \phi(w)}%
    {\label{eq:RegularizedNLP}}{}%
    \addConstraint{x_0}{=  \bar x_0}%
    \addConstraint{0}{      = G_\mathrm{IRK}(x_k, x_{k+1}, u_k, z_k), \quad & \forall k \in \mathcal{K} }%
    \addConstraint{0}{\leq h(x_k, u_k, z_k),  \quad & \forall k \in \mathcal{K}}%
\end{mini!}%
For the case of $p,q \rightarrow\infty$, the regularization cost \eqref{eq:pqregularization} implicitly implements the constraint:
\begin{align}
    -E_\mathrm{max} \leq \hat{E}(w) \leq E_\mathrm{max} \jochem{\ .}
\end{align}
In order to avoid issues with infeasibility, finite values of $p$ and $q$ are recommended in practice.

\begin{figure}
    \centering
    \includegraphics[width=1\linewidth, trim=0.2cm 0.4cm 0.2cm 0, clip]{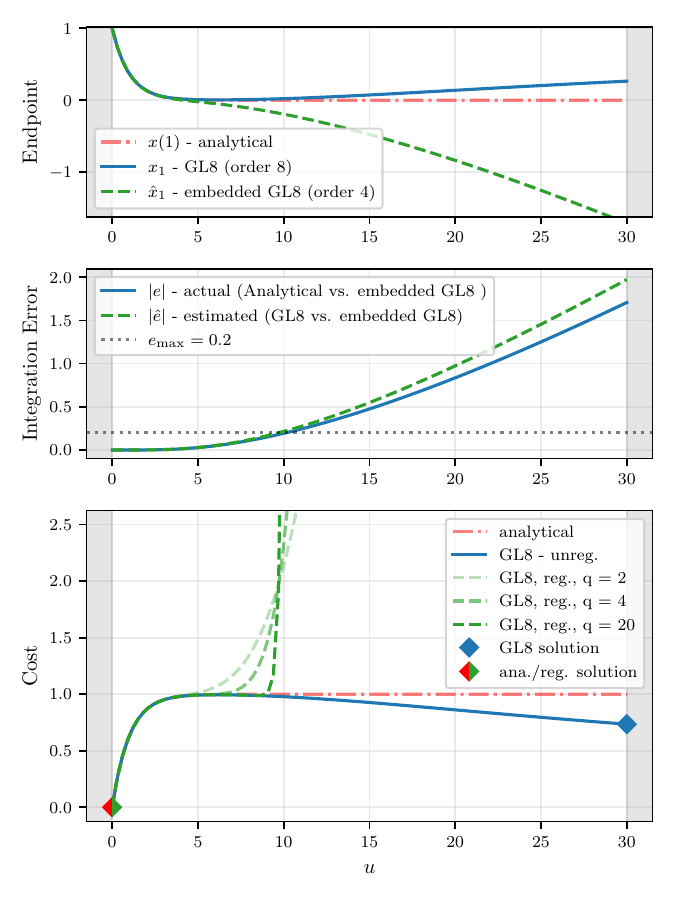}
    \caption{Minimal example, cont.: Endpoint value, integration error and objective over the feasible set. The difference between the endpoints of the embedded and original integrator is used to regularize the objective.}
    \label{fig:SimpleExample_RegOrConstr}
\end{figure}

\medskip\noindent\textbf{(Minimal Example, cont.)} To estimate the integration error in the single integration interval in the fashion of Eq. \eqref{eq:embIntErrApprox}, we add a (lower-order) embedded integrator to the GL8 scheme in the discretized problem \eqref{eq:simpleExample_OG}. Fig. \ref{fig:SimpleExample_RegOrConstr} shows for different control values the integration endpoint $x_1$ for both the GL8 and the embedded integrator (top) and shows the corresponding estimation of the local integration error compared to the real local integration error (middle).
Fig.~\ref{fig:SimpleExample_RegOrConstr} (bottom) shows the NLP cost function including the regularization term \eqref{eq:pqregularization} with $p = 2$, $\omega = 1$ and $e_{\max} = 0.2$ and different values of $q$.
For this value of $e_{\max}$, the solver converges to the correct solution, independent of the initialization.

\subsection{Discussion}

By adding the regularization discussed above we modify the cost function of the discretized NLP to penalizing trajectories with a large integration error. 
This might result in less optimal but higher accuracy trajectories.
The proposed method introduces only a single tuning factor (the maximum estimated relative error), which is intuitively to choose.

This `numerical robustification' comes at a numerical cost: The embedded integrator requires at least one extra stage where the dynamics are evaluated.
If the embedded integrator pair is chosen wisely, for example by adding an extra stage $c_0 = 0$ to a Radau2A-3 scheme, the stage value $z_{k,0}$ is simply the shooting node $x_k$. If now also the control parametrization is continuous, no extra dynamics evaluation is required for this stage since the evaluation from the previous interval can be used. 

Also we would like to point out that the proposed regularization is based on the error estimate \eqref{eq:embIntErrApprox} which is an estimate for local error of the lower order embedded integrator. Thus regularizing this estimate does not guarantee a reduction of the local error of the higher-order integrator, whose endpoint we use to continue the integration.

The expression of $e_k$ is then a possibly nonlinear function of the stage variables and the control, which makes the resulting NLP more nonlinear.
Although the goal of the regularization is to rule out spurious minima, the problem might still appear the if the integration error of both integrators aligns.




\section{Example: Hang Glider Problem}
As a more realistic example, consider the `Hang Glider' example from \cite{Betts2010}, originally posed by \cite{Bulirsch1993},
of controlling a hang glider that flies over a thermal updraft.
For the problem, the glider systems state contains the horizontal and vertical position and velocity, the control is the aerodynamic lift coefficient of the glider:
\begin{align}
     x = (p_x, p_y, v_x, v_y)^\top \in \R^4, &&  u = C_L \in \R
\end{align}
For the details about the stiff point-mass state dynamics $f: \R^{4} \times \R \rightarrow \R^4 $, the details on the updraft formulation and the chosen parameter values, we refer the reader to \cite{Betts2010}.

We want to find the optimal control strategy and final time which maximizes the horizontal distance that the glider flies from a fixed initial position and velocity $\bar{x}_0 = [\bar{p}_{x,\mathrm{0}},\,\bar{p}_{y,\mathrm{0}},\,\bar{v}_{x,\mathrm{0}},\,\bar{v}_{y,\mathrm{0}}]^\top =
    \SI{0}{\meter},\,
    \SI{1000}{\meter},\,
    \SI{13.277}{\meter\per\second},\,
    \SI{-1.2875}{\meter\per\second}]^\top $ to a target vertical height and velocity. 
    The target endpoint is originally posed as a terminal constraint, here we express it with the terminal cost
\begin{align}
    E(x(t_\mathrm{f})) = \left\|\frac{\tilde{\Omega}^{-1}}{10^{-2}}\left(\begin{bmatrix}
        p_y(t_\mathrm{f}) \\
        v_x(t_\mathrm{f}) \\
        v_y(t_\mathrm{f})\end{bmatrix} - \begin{bmatrix}
        \SI{900}{\meter} \\
        \bar{v}_{x,\mathrm{0}} \\
        \bar{v}_{y,\mathrm{0}}
    \end{bmatrix}\right)\right\|^2_2
\end{align}
where $\tilde{\Omega} = \mathrm{diag}(\SI{100}{\meter},\SI{10}{\meter\per\second},\SI{5}{\meter\per\second})$ normalizes the state values with respect to their expected range.
This later allows us to compare exact and inexact solution by their objective value.
The continuous time OCP is then given by
\begin{mini!}
    {x(\cdot), u(\cdot), t_\mathrm{f}}{-\frac{1}{\omega_{p,x}}p_x(t_\mathrm{f}) + E(x(t_\mathrm{f}))\label{ocp:GliderOCP_objective}}
    {\label{ocp:GliderOCP}}{}
    \addConstraint{0}{= x(0) - \bar x_0}
    \addConstraint{\dot x(t)}{= f(x(t), u(t)), \qquad&&\forall t\in [0,t_\mathrm{f}]}
    \addConstraint{0}{\leq u(t) \leq 1.5, \qquad&&\forall t\in [0,t_\mathrm{f}]} 
\end{mini!}
where later we will use $J: \R^{n_x} \rightarrow \R$ for the objective \eqref{ocp:GliderOCP_objective} which only depends on the endpoint of the state trajectory, the weight $\omega_{p,x} = \SI{1000}{\meter}$ normalizes the objective.
To solve the problem, we use the following discretization: We divide the horizon $t \in [0,t_\mathrm{f}]$ into $N=30$ intervals in each of which we use a constant control $u_k$.
For this choice of control parametrization, the optimal achievable objective is $p_x^* = \SI{1247}{\meter}$, a value slightly smaller than the one reported by \cite{Betts2010} who used repeated refinements of the control and state mesh.

In every interval we use a single step of either a Heun-Euler or the Radau2A-3 scheme with an embedded integrator.
The free end time is implemented by scaling the dynamics as $\tilde{f}(x,u,t_\mathrm{f}) = t_\mathrm{f} f(x,u)$, considering the numerical timescale $\tau \in [0,1]$ and adding $t_\mathrm{f}$ to the decision variables of the NLP.
The variables are initialized according to \cite{Betts2010}. We use IPOPT\cite{IPOPT} equipped with the linear solver MA27\cite{HSLMA27} to solve the NLPs formulated using CasADi\cite{Andersson2019} in Python 3.8.


We regularize the estimated integration error with the additional cost term \eqref{eq:pqregularization}, where we choose $p,q = 2$ and $E_\mathrm{max}(e_\mathrm{max})$ in accordance to equation \eqref{eq:Emax_singleTuning} with state weight vector $\omega = [\SI{1000}{\meter},\SI{100}{\meter},\SI{10}{\meter\per\second},\SI{5}{\meter\per\second}]^\top$, such that we only have the single tuning parameter $e_\mathrm{max}$.

To quantify the performance of a found solution to the NLP, we perform a high-accuracy simulation denoted by $x_\mathrm{sim}(t)$ of the found control strategy and compare the simulation values on the grid points $t_1, \dots, t_N$ with the computed shooting variables, to compute the average simulation error:
\begin{align}
    E_\mathrm{sim} := \frac{1}{N} \sum_{k \in \mathcal{K}}^{ }\left\| \Omega^{-1} \left(x_\mathrm{sim}(t_{k+1}) - x_{k+1})\right) \right\|_2
\end{align}
where $\Omega = \mathrm{diag}(\omega)$ is formed using the weighting vector from above.
For the simulation we can also evaluate the objective $J_\mathrm{sim} := J(x_\mathrm{sim}(t_\mathrm{f}))$.

\begin{figure}
    \includegraphics[width=\linewidth, trim=0.2cm 0.4cm 0 0, clip]{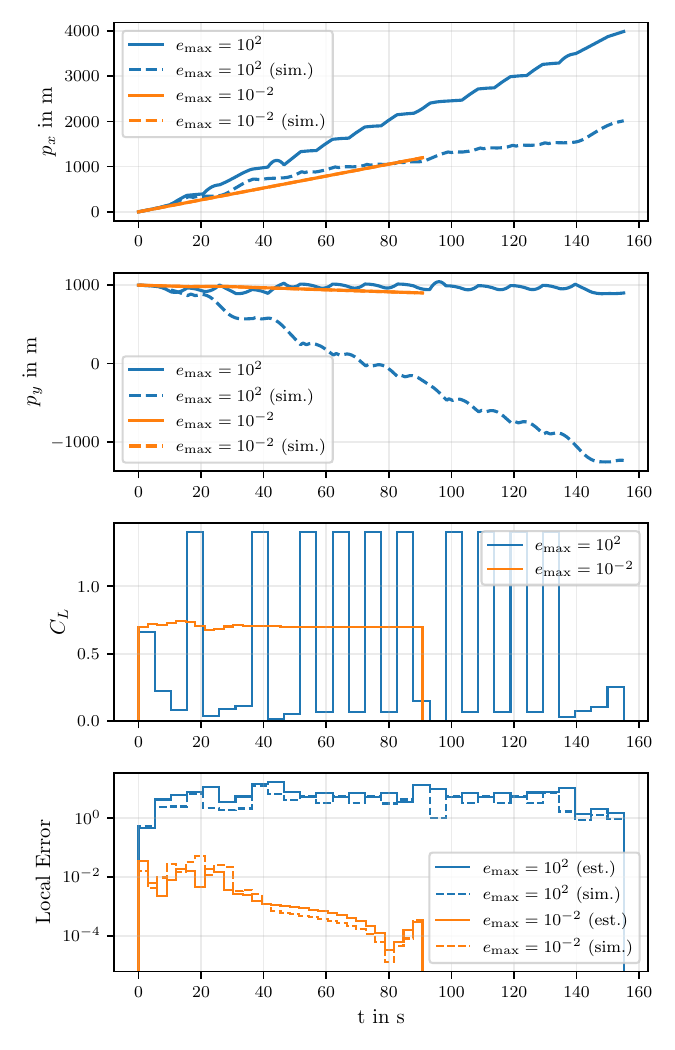}
    \caption{Solution of the discretized hang glider problem using a Heun-Euler integration scheme, unregularized  $e_\mathrm{rel,max} = 10^{2}$ (blue) and regularized $e_\mathrm{rel,max} = 10^{-2}$ (orange). To avoid a large integration error, a less aggressive control strategy is optimal.}
    \label{fig:Glider_HE_comparison}
\end{figure}

\begin{figure}
    \includegraphics[width=\linewidth, trim=0.2cm 0.4cm 0 0, clip]{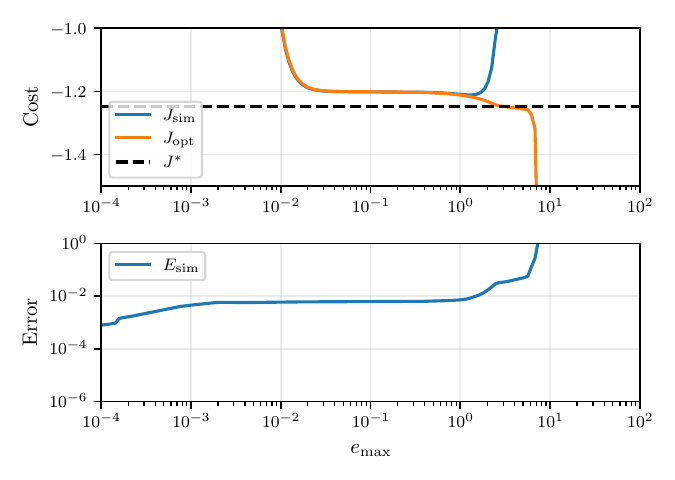}
    \caption{Sweep of the regularization parameter $e_\mathrm{rel,max} \in [10^{-4}, 10^{2}]$, for the embedded explicit Heun-Euler integration scheme with $N=30$ intervals. We can make a very cheap integration method work.}
    \label{fig:Glider_HE_sweep}
\end{figure}

Figure \ref{fig:Glider_HE_comparison} compares the unregularized ($e_\mathrm{max} = 10^{2}$) and the regularized ($e_\mathrm{max} = 10^{-2}$) solutions to the discretized problem using Heun-Euler integration. 
In the unregularized case, due to the `bad' discretization, it is optimal to aggressively operate the control (third plot), which leads to a high integration error (last plot), which minimizes the objective (first plot).
This prediction is not accurate, the comparison simulation shows a completely different behavior (dashed line).
By choosing a maximum estimated integration error of $e_\mathrm{max} = 10^{-2}$, we converge to a solution that matches the posterior simulation quite well.
This achieved by less aggressive control maneuvers, which reduces the estimated error of the embedded Euler method (last plot).
In the last plot, the dashed line represent the actual local error of the Euler method in each interval.

Figure \ref{fig:Glider_HE_sweep} shows a sweep of the regularization parameter $e_\mathrm{max} \in [10^{-4}, 10^{2}]$ to compare optimality and integration error $E_\mathrm{sim}$ of the solution.
For $e_\mathrm{max} \leq 1$, the OCP converges to a suboptimal but accurate solutions. For larger values, the solutions are no longer accurate.

\begin{figure}
    \includegraphics[width=\linewidth, trim=0.2cm 0.4cm 0 0, clip]{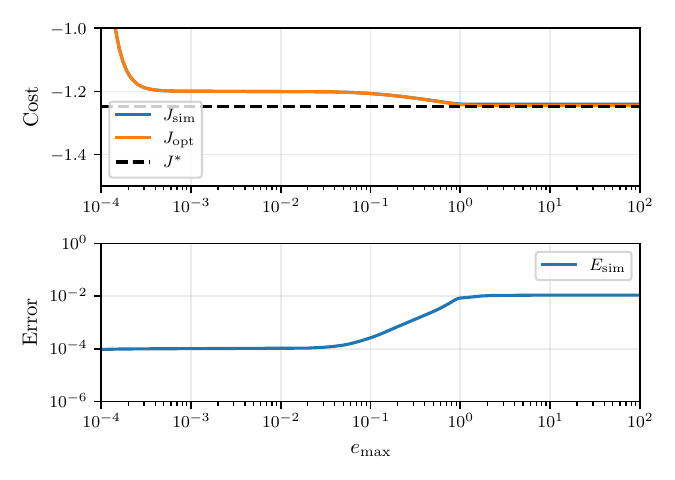}
    \caption{Sweep of the regularization parameter $e_\mathrm{rel,max} \in [10^{-4}, 10^{2}]$, for the embedded Radau2A-3 integration scheme with $N=30$ intervals. We can trade off optimality with accuracy.}
    \label{fig:Glider_RD_sweep}
\end{figure}

\begin{figure}
    \includegraphics[width=\linewidth, trim=0.2cm 0.4cm 0 0, clip]{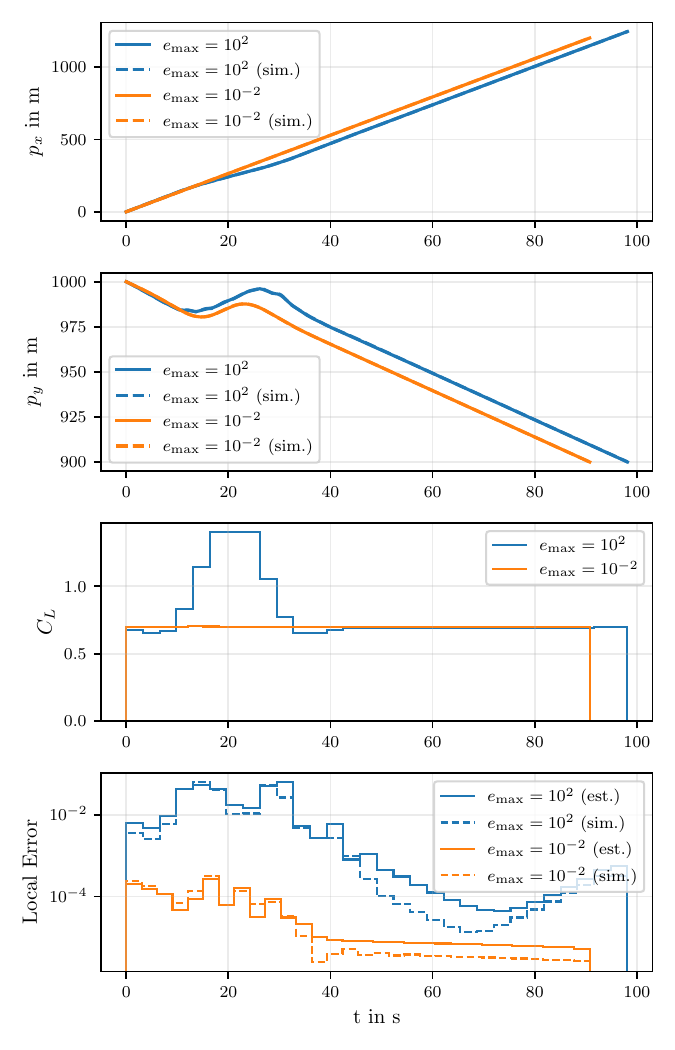}
    \caption{Solution of the regularized discretized problem using a Radau2A-3 integration scheme, with $e_\mathrm{rel,max} = 10^{2}$ (blue) and $e_\mathrm{rel,max} = 10^{-2}$ (orange). To avoid a large integration error, a less aggressive control strategy is optimal.}
    \label{fig:Glider_RD_comparison}
\end{figure}

In comparison with the very cheap Heun-Euler integration scheme, the Radau2A-3 scheme delivers, even for the unregularized case, consistent results with a maximum error score of $10^{-2}$, compare Figure \ref{fig:Glider_RD_sweep}.
By regularization of the approximated integration error, we can further reduce this simulation error for the cost of optimality. 
Again, this is achieved by choosing a much less aggressive control strategy as can be seen in Figure \ref{fig:Glider_RD_comparison}.

Comparing the efficiency and the optimality for a Radau2A-3 discretization without regularization with a Heun-Euler discretization with the best possible regularization $e_\mathrm{max} = 0.1$, we observe a speedup in time-per-iteration from $\SI{1.2}{\milli\second}$ to $\SI{0.4}{\milli\second}$ at a cost increase of $-1.24$ to $-1.20$ ($\SI{3}{\percent}$), compare Figures \ref{fig:Glider_HE_sweep} and \ref{fig:Glider_RD_sweep}.


\section{Conclusion}

We presented how embedded Runge-Kutta methods can be used to estimate and regularize the local integration error of a discretized optimal control problem. 
On an example problem, we showed that using the regularization, we can trade off optimality and accuracy for a given discretization, and also enable the use of cheap, low-order integration schemes that are otherwise unsuitable.

In future work we will investigate better error estimates, the computational efficiency gains, and compare other methods of integration error estimate.


\section*{Acknowledgements}
The authors would like to thank Amon Lahr, Armin Nurkanovic and Per Rutquist for the fruitful discussions.

\bibliography{bibliography/sources.bib}


\end{document}